\documentclass[12pt,a4paper]{amsart}

\usepackage{amsmath, amssymb, amsthm}
\theoremstyle{plain}
\newtheorem{thm}{Theorem}[section]
\newtheorem{lem}[thm]{Lemma}
\newtheorem{prop}[thm]{Proposition}

\newtheorem{THM}{Theorem}

\theoremstyle{remark}

\usepackage{txfonts}
\usepackage{mathbbol}
\usepackage{ae}

\usepackage{color}
\usepackage{graphicx}
\usepackage{url}
\usepackage[abs]{overpic}
\usepackage{subfigure}
\usepackage{booktabs}
\usepackage{paralist}
\usepackage[a4paper]{anysize}

\newcommand\interior[1]{\operatorname{int}{#1}}

\renewcommand{\phi}{\varphi}

\begin{document}

\title[Genocchi Numbers and Simplicial Balls]{Genocchi Numbers and $f$-Vectors of Simplicial Balls}
\author[S. Herrmann]{Sven Herrmann}
\address{Sven Herrmann, Fachbereich Mathematik, AG~7, Technische 
Universität Darmstadt, 64289 Darmstadt, Germany}
\email{sherrmann@mathematik.tu-darmstadt.de}
\date{\today}

\begin{abstract}
The \emph{Genocchi numbers} are the coefficients of the generating function $\frac{2t}{e^t+1}$. In this note we will give an equation for simplicial balls which involves this numbers. It relates the number of faces in the interior of the ball with the number of faces in the boundary of the ball. This is a variation of similar equations given in \cite{MR2070631} and \cite{HJ}.
\end{abstract}

\maketitle

\section{Introduction}
\noindent
The aim of this note is to investigate $f$-vectors of simplicial balls, especially the relations between interior and boundary faces. For a simplicial ball $B$ we denote by $f_i(B)$ the number of $i$-dimensional faces. The boundary $\partial B$ of $B$ is a simplicial sphere with face numbers 
$f_i(\partial B)$. We also define $f_i(\interior B):=f_i(B)-f_i(\partial B)$ 
although the interior $\interior B$ of $B$ is not a polyhedral complex.

For simplicial spheres the \emph{Dehn-Sommerville relations} give non-trivial dependences for the $f$-vector. McMullen and Walkup gave in \cite{MR0298557} a similar equation for simplicial balls which relates the $f$-vector of the boundary $\partial B$ with the $f$-vector of $B$ (see also Billera and Björner \cite{MR1730171} and Klain \cite{Klain}). The most important examples of simplicial balls are triangulated polytopes. Special interest on connections between the $f$-vector of the boundary and the interior occurs for example when the polytope is simplicial and so the $f$-vector of the boundary is fixed (see McMullen \cite{MR2070631}), but also in other cases, e.g. for triangulations without interior edges of low dimension \cite{HJ}.

The $f$-vector of the interior of a simplicial ball occurs in the study of the combinatorics of certain unbounded polyhedra \cite[Section 2]{HJ}: It is dual to the $f$-vector of the polyhedral complex of their bounded faces. A very important special case is the tight span of a finite metric space \cite{MR753872,HJ}.

We will give a relation between the $f$-vector of the boundary and the interior of a simplicial ball directly in terms of the $f$-vector of the interior. The most interesting point about our equation is the occurrence of the Genocchi numbers $G_{2n}$:

\begin{THM}
  For each simplicial $(n-1)$-ball~$B$ and $n-k$ even we have
  \[
  f_k(\interior{B})\ =\ \sum_{i=1}^{\lfloor\frac{n-k}2\rfloor}\frac{G_{2i}}{2i}\left(
    \binom{k+2i-1}{k+1}f_{k+2i-2}(\partial B)-\binom{k+2i}{k+1} f_{k+2i-1}(\interior B)\right) \quad .
  \]
\end{THM}

The proof of this theorem is given in Section 3, while Section 2 collects basic facts about the Genocchi numbers.

\section{Genocchi Numbers}

\noindent
The \emph{Genocchi numbers} (Sequence A001469 in Sloane \cite{OEIS}) are 
defined by means of the generating function

\begin{align}\label{gn}
\frac{2t}{e^t+1}=:\sum_{n=0}^\infty G_n \frac{t^n}{n!}=t+\sum_{n=1}^\infty G_{2n} \frac{t^{2n}}{(2n)!}\quad.
\end{align}
Their investigation goes back to Leonhard Euler while the term Genocchi numbers (after the Italian mathematician Angelo
Genocchi (1817--1889)) was presumably introduced by Edouard Lucas in \cite{Lucas}. Intensive studies on these numbers
were done by Eric T. Bell in the 1920s \cite{MR1501336,MR1501490}.

The Genocchi numbers are directly related to the much more famous \emph{Bernoulli numbers} $B_n$ via
$G_n=2(1-2^{2n})B_n$. There is also a combinatorial interpretation which is due to Dominique Dumont \cite{MR0337643}:
The absolute value of $G_{2n+2}$ equals the number of permutations $\tau$ of $\{1,2,\dots,2n\}$ such that $\tau(i)>i$ 
if and only if $i$ is odd. There is a continuing interest in 
the Genocchi numbers and especially their generalizations (see e.g. \cite{MR2020979,MR2110777,MR2206473}).

There are a lot of possibilities of computing the values of the Genocchi numbers (see for example \cite{MR2110777,MR1771988,TT}, and also the entry in \cite{OEIS} for further references). Here we need two other formulae which are obtained by differentiating \eqref{gn} and comparing the coefficients of even and odd powers of $t$ respectively.

\begin{prop}\label{prop:genocchi}
  For the Genocchi numbers we have the recursion formulae
  \[
  G_{2n}=-n-\frac1 2 \sum_{k=1}^{n-1}\binom{2n}{2k}G_{2k}\quad
  \text{and}\quad
  G_{2n}=-1-\sum_{k=1}^{n-1}\binom{2n}{2k-1}\frac{G_{2k}}{2k}\quad.
\]
\end{prop}

Additionally, we will use the following technical lemma:
\begin{lem}\label{lem:genocchi}
  For $n \geq 2$ we have
  \[
  \frac{2n-1}{n}G_{2n}=-\frac12\sum_{k=1}^{n-1}\binom{2n-1}{2k-1}\frac{2k-1}k G_{2k}\quad.
  \]
\end{lem}
\begin{proof}
  First we remark that the second part of Proposition \ref{prop:genocchi} yields
  \begin{equation}\label{eq:gen}
    1=-\sum_{k=1}^{n-1}\binom{2n-2}{2k-1}\frac{G_{2k}}{2k}
  \end{equation}
  for all $n\geq2$. Then we use the first part of Proposition \ref{prop:genocchi} and compute
  \begin{align*}
    \frac{2n-1}{n}G_{2n}\ \stackrel{\text{\ref{prop:genocchi}}}{=}\ & -(2n-1)-\frac12\sum_{k=1}^{n-1}\binom{2n}{2k}\frac{2n-1}nG_{2k} \\
    \stackrel{\text{\eqref{eq:gen}}}{=}\ &(2n-1)\sum_{k=1}^{n-1}\binom{2n-2}{2k-1}\frac{G_{2k}}{2k}-\frac12\sum_{k=1}^{n-1}\binom{2n}{2k}\frac{2n-1}nG_{2k}\\
    =\ &-\frac12\sum_{k=1}^{n-1}\binom{2n-1}{2k-1}\frac{2k-1}k\left(\frac{2n-1}{2k-1}-\frac{2n-2k}{2k-1}\right)G_{2k}\\
    =\ &-\frac12\sum_{k=1}^{n-1}\binom{2n-1}{2k-1}\frac{2k-1}k G_{2k}\quad.
  \end{align*}
\end{proof}

\section{Proof of Theorem 1}
\noindent
To prove our main theorem we use a variation of the Dehn-Sommerville relations by Klain \cite{Klain} for arbitrary triangulated manifolds with boundary. The following is a reformulation of \cite[Corollary 1.4]{Klain}. (Note that there is a typing error in Klain's version: the $-1$ in the exponent of $-1$ should be left out.)

\begin{prop}\label{prop:DS-f-vector}
  For each simplicial $(n-1)$-ball~$B$ and $n-k$ even we have
  \[
  f_k(\interior{B})\ =\ -\frac{f_k(\partial B)}2 - \sum_{i=1}^{n-k-1}\frac{(-1)^{n+k+i}}{2}\binom{k+1+i}{k+1}f_{k+i}(\interior{B})\quad.
  \]
\end{prop}

Now we can give the proof of our main result.

\begin{proof}[Proof of Theorem 1]
  The claim is trivial for $k=n$ so we may assume that it holds for $k'>k$. For simplicity we also assume that $n$ (and
  thus $k$) is even.  The case $n$ odd is analogous.  Thus we compute

  \begin{align*}
    f_k(\interior{B})\ \stackrel{\ref{prop:DS-f-vector}}=\ &-\frac{f_k(\partial B)}2-\sum_{i=1}^{n-k-1}\frac{(-1)^{n+k+i}}{2}\binom{k+1+i}{k+1}f_{k+i}(\interior{B})\\
    \ =\ & -\frac{f_k(\partial B)}2+\frac12\sum_{i=1}^{\frac{n-k}{2}}\binom{k+2i}{k+1}f_{k+2i-1}(\interior{B})-\frac12\sum_{i=1}^{\frac{n-k}{2}-1}\binom{k+2i+1}{k+1}f_{k+2i}(\interior{B})\\
    \ =\ & -\frac{f_k(\partial B)}2+\frac12\sum_{i=1}^{\frac{n-k}{2}}\binom{k+2i}{k+1}f_{k+2i-1}(\interior{B})-\frac{1}{2}\sum_{i=1}^{\frac{n-k}{2}-1}\binom{k+2i+1}{k+1}\\
    & \cdot \sum_{j=1}^{\frac{n-k}{2}-i}\frac{G_{2j}}{2j}\left( \binom{k+2i+2j-1}{k+2i+1}f_{k+2i+2j-2}(\partial B)-\binom{k+2i+2j}{k+2i+1} f_{k+2i+2j-1}(\interior B)\right)\\ 
    \ \overset{\text{($*$)}}=\ &-\frac{f_k(\partial B)}2+\frac12\sum_{i=1}^{\frac{n-k}{2}}\binom{k+2i}{k+1}f_{k+2i-1}(\interior{B})-\frac{1}{2}\sum_{j=2}^{\frac{n-k}{2}}\sum_{i=1}^{j-1}\binom{k+2j-2i+1}{k+1}\\
    &\quad \cdot\frac{G_{2i}}{2i}\left( \binom{k+2j-1}{2i-2}f_{k+2j-2}(\partial B)-\binom{k+2j}{2i-1} f_{k+2j-1}(\interior B)\right)\\
    =\ &-\frac{f_k(\partial B)}2-\sum_{j=2}^{\frac{n-k}{2}}\binom{k+2j-1}{k+1}f_{k+2j-2}(\partial B)\left(\frac12\sum_{i=1}^{j-1}\binom{2j-1}{2i-1}\frac{2i-1}{(2j-1)i}G_{2i}\right)\\
    &\quad +\sum_{j=1}^{\frac{n-k}{2}}\frac1{2j}\binom{k+2j}{k+1} f_{k+2i-1}(\interior B)\left(j+\frac12\sum_{i=1}^{j-1}\binom{2j}{2i}G_{2i}\right)	\\
    \ \overset{\substack{\ref{lem:genocchi},\\\ref{prop:genocchi}}}=\ &\sum_{j=1}^{\frac{n-k}2}\frac{G_{2j}}{2j}\left( \binom{k+2j-1}{2j-2}f_{k+2j-2}(\partial B)-\binom{k+2j}{2j-1} f_{k+2j-1}(\interior B)\right)\quad.
  \end{align*}
  Where in  ($*$) we used the computation 
  \[
  \sum_{i=1}^{\frac{n-k}{2}-1}\sum_{j=1}^{\frac{n-k}{2}-i}\alpha(i,j,k)=\sum_{i=1}^{\frac{n-k}{2}-1}\sum_{j=i+1}^{\frac{n-k}{2}}\alpha(i,j-i,k)=\sum_{j=2}^{\frac{n-k}{2}}\sum_{i=1}^{j-1}\alpha(i,j-i,k)=\sum_{j=2}^{\frac{n-k}{2}}\sum_{i=1}^{j-1}\alpha(j-i,i,k)\quad .
  \]
\end{proof}
If we now  consider a simplicial $(n-1)$-ball  $B$ without any interior faces of dimension up to $e$ we derive
from this theorem the equations
\[
  \sum_{i=1}^{\lfloor\frac{n-k}2\rfloor}\frac{G_{2i}}{2i}\binom{k+2i-1}{k+1}f_{k+2i-2}(\partial B)
  \ =\ \sum_{i=1}^{\lfloor\frac{n-k}2\rfloor}\frac{G_{2i}}{2i}\binom{k+2i}{k+1} f_{k+2i-1}(\interior B)
\]
for $n-k$ even and $k\leq e$. So we have $\lfloor (e+2)/2\rfloor$ ($\lfloor (e+3)/2\rfloor$) equations for the
$\lfloor (n-1-e)/2\rfloor$ ($\lfloor (n-e)/2\rfloor$) unknowns $f_k(\interior{B})$ with $n-k$ odd for $n$ even and odd
respectively. This is an $f$-vector analog of Proposition 3.4 in \cite{HJ}.

\bibliographystyle{abbrv}
\bibliography{main}

\end{document}